\documentclass[11pt,reqno]{amsart}

\usepackage[latin1]{inputenc}
\usepackage[left=3.5cm,right=3.5cm,top=3.5cm,bottom=3.5cm]{geometry}
\usepackage{amssymb}
\usepackage{graphicx}
\usepackage{amscd}
\usepackage[hidelinks]{hyperref}
\usepackage{color}
\usepackage{float}
\usepackage{graphics,amsmath,amssymb}
\usepackage{amsthm}
\usepackage{amsfonts}
\usepackage{latexsym}
\usepackage{epsf}
\usepackage{enumerate}
\usepackage{xifthen}
\usepackage{mathrsfs}
\usepackage{dsfont}
\usepackage{makecell}
\usepackage[FIGTOPCAP]{subfigure}
\usepackage{amsmath}
\allowdisplaybreaks[4]
\usepackage{listings}
\usepackage{etoolbox}

 \newtheoremstyle{mytheorem}
 {3pt}
 {3pt}
 {\slshape}
 {}
 {\bfseries}
 {.}
 { }
 {}

\numberwithin{equation}{section}

\theoremstyle{mytheorem}
\newtheorem{theorem}{Theorem}[section]

\newtheorem{lemma}[theorem]{Lemma}

\theoremstyle{definition}

\newtheorem{remark}{Remark}[section]

\newcommand{\doi}[1]{\href{http://dx.doi.org/#1}{doi: #1}}

\newcommand{\Keywords}[1]{\ifthenelse{\isempty{#1}}{}{\smallskip \smallskip \noindent \textbf{Keywords}. #1}}
\newcommand{\MSC}[2][2010]{\ifthenelse{\isempty{#2}}{}{\smallskip \smallskip \noindent \textbf{#1MSC}. #2}}

\makeatletter
\patchcmd{\@settitle}{\uppercasenonmath\@title}{\LARGE}{}{}
\patchcmd{\@setauthors}{\MakeUppercase}{\large}{}{}
\makeatother

\title[Formulas for partition $k$-tuples with $t$-cores]{Formulas for partition $\boldsymbol{k}$-tuples with $\boldsymbol{t}$-cores} 

\author[S. Chern]{Shane Chern}
\address{School of Mathematical Sciences, Zhejiang University, Hangzhou, 310027, China}
\email{\href{mailto:shanechern@zju.edu.cn}{shanechern@zju.edu.cn}; \href{mailto:chenxiaohang92@gmail.com}{chenxiaohang92@gmail.com}}

\date{}

\begin{document}

{\footnotesize\noindent \textit{J. Math. Anal. Appl.} \textbf{437} (2016), no. 2, 841--852.\\
\doi{10.1016/j.jmaa.2016.01.040}}

\bigskip \bigskip

\maketitle

\thispagestyle{empty}

\begin{abstract}  
Let $A_{t,k}(n)$ denote the number of partition $k$-tuples of $n$ where each partition is $t$-core. In this paper, we establish formulas of $A_{t,k}(n)$ for some values of $t$ and $k$ by employing the method of modular forms, which extends Wang's result for $t=3$ and $k=2,3$.

\Keywords{Partition, $k$-tuple, $t$-core, modular form.}

\MSC{Primary 11P84; Secondary 11F11, 11M36, 05A17.}
\end{abstract}

\section{Introduction}

Let $A_t(n)$ denote the number of $t$-core partitions of $n$, that is, the number of partitions of $n$ with no hook numbers being multiples of $t$. Garvan, Kim, and Stanton \cite[Eq. (2.1)]{GKS1990} showed that the generating function of $A_t(n)$ is given by
\begin{equation}\label{eq:1.1}
\sum_{n\ge 0}A_t(n)q^n=\frac{(q^t;q^t)_\infty^t}{(q;q)_\infty},
\end{equation}
where, as usual, we denote
$$(a;q)_\infty:=\prod_{n\ge 0}(1-aq^n),$$
and
$$(a;q)_n:=\frac{(a;q)_\infty}{(aq^n;q)_\infty}\quad (-\infty<n<\infty).$$
For convenience, we also write
$$(a_1,a_2,\ldots,a_n;q)_\infty:=(a_1;q)_\infty(a_2;q)_\infty\cdots (a_n;q)_\infty.$$

We say $(\lambda_1,\ldots,\lambda_k)$ is a partition $k$-tuple of $n$ if the sum of all the parts equals $n$. For example, $(\{1,1\},\{1\})$ is a partition pair of $3$. Furthermore, a partition $k$-tuple $(\lambda_1,\ldots,\lambda_k)$ of $n$ with $t$-cores means that each $\lambda_i$ is $t$-core. Let $A_{t,k}(n)$ denote the number of partition $k$-tuples of $n$ with $t$-cores. From \eqref{eq:1.1}, we readily obtain the generating function of $A_{t,k}(n)$, that is,
\begin{equation}\label{eq:1.2}
\sum_{n\ge 0}A_{t,k}(n)q^n=\left(\sum_{n\ge 0}A_t(n)q^n\right)^k=\frac{(q^t;q^t)_\infty^{kt}}{(q;q)_\infty^k}.
\end{equation}
Here we write $A_{t,1}(n)=A_t(n)$.

Many authors have studied the number of partitions and partition pairs with $t$-cores and obtained sets of Ramanujan-like congruences (see, e.g., \cite{BN2014, BN2015, Boy2002, Chen2013, Dai2015, HS1996, Lin2014, Yao2015}). More recently, Wang \cite{Wang2016} established formulas of $A_{3,2}(n)$ and $A_{3,3}(n)$ with the help of Ramanujan's $_1\psi_1$ formula and Bailey's $_6\psi_6$ formula. Let $\sigma_{1}(n)=\sum_{d\mid n}d$ and $\sigma_{2,\chi_3}^{*}(n)=\sum_{d\mid n}\chi_3(\frac{n}{d})d^2$ (here $\chi_3(n)=(n|3)$ denotes the Legendre symbol). He proved that for any integer $n\ge 0$,
\begin{equation}\label{eq:wang1}
A_{3,2}(n)=\frac{1}{3}\sigma_{1}(3n+2),
\end{equation}
and
\begin{equation}\label{eq:wang2}
A_{3,3}(n)=\sigma_{2,\chi_3}^{*}(n+1).
\end{equation}

However, we notice that Wang's method expires for $k\ge 4$ due to a lack of corresponding $_s\psi_s$ formulas. Recall that for the case $k=1$, Granville and Ono \cite{GO1996} gave the formula
\begin{equation}\label{eq:GO}
A_{3,1}(n)=\sigma_{0,\chi_3}(3n+1),
\end{equation}
where $\sigma_{0,\chi_3}(n)=\sum_{d\mid n}\chi_3(d)$, using the tools of modular forms. We also notice that in \cite{ORW1995}, Ono, Robins, and Wahl found formulas of the number of representations of $n$ as sums of $k$ triangular numbers (denoted by $\delta_k(n)$) for $k=2$, $4$, $6$, $8$, $10$, $12$, and $24$ by applying the same method. It is known by Jacobi's identity that
$$\sum_{n\ge 0}q^{T_n}=1+q+q^3+q^6+q^{10}+\cdots=\frac{(q^2;q^2)_\infty^2}{(q;q)_\infty}.$$
We immediately see that their $\delta_k(n)$ is our $A_{2,k}(n)$. It is therefore natural to expect the method of modular forms will play a role in finding formulas of $A_{t,k}(n)$ for some other values of $t$ and $k$.

In Section \ref{sec:2}, we will give a brief introduction to modular forms. In Section \ref{sec:3}, 
we will discuss $A_{3,k}(n)$ and derive two new formulas for $k=4$ and $6$. We should say that these formulas, which involve Fourier coefficients of some $\eta$-products, are not explicit. However, they can be regarded as analogues to $\delta_k(n)$ in \cite{ORW1995}. In the following sections, we will present some formulas for $t\ge 4$. Here the formula of $A_{5,1}(n)$ is explicit.

\section{A brief introduction to modular forms}\label{sec:2}

Let $\gamma=\begin{pmatrix}a & b\\c & d\end{pmatrix}\in SL_2(\mathbb{Z})$. We first define the following congruence subgroups of level $N$:
\begin{enumerate}
\item $$\Gamma_0(N)=\left\{\gamma\in SL_2(\mathbb{Z}): \gamma\equiv \begin{pmatrix}* & *\\0 & *\end{pmatrix}\pmod{N}\right\}.$$
\item $$\Gamma_1(N)=\left\{\gamma\in SL_2(\mathbb{Z}): \gamma\equiv \begin{pmatrix}1 & *\\0 & 1\end{pmatrix}\pmod{N}\right\}.$$
\item $$\Gamma(N)=\left\{\gamma\in SL_2(\mathbb{Z}): \gamma\equiv \begin{pmatrix}1 & 0\\0 & 1\end{pmatrix}\pmod{N}\right\}.$$
\end{enumerate}
Here ``$*$'' means ``unspecified.''

Let $A\in SL_2(\mathbb{Z})$ act on the complex upper half plane $\mathcal{H}$ by the linear fractional transformation
$$A\tau=\frac{a\tau+b}{c\tau+d}.$$
Let $\chi$ be a Dirichlet character mod $N$ and $k\in\mathbb{Z}^+$ satisfying $\chi(-1)=(-1)^k$. Let $f(\tau)$ be a holomorphic function on $\mathcal{H}$ such that
$$f(A\tau)=\chi(d)(c\tau+d)^kf(\tau)$$
for all $A\in\Gamma_0(n)$ and all $\tau\in\mathcal{H}$. We call $f(\tau)$ a modular form of weight $k$ and nebentypus $\chi$ on $\Gamma_0(N)$. Furthermore, we say $f(\tau)$ is a holomorphic modular form (resp. cusp form) if $f(\tau)$ is holomorphic (resp. vanishes) at the cusps of $\Gamma_0(N)$. It is known that the holomorphic modular forms (resp. cusp forms) of weight $k$ and nebentypus $\chi$ form finite dimensional $\mathbb{C}$-vector spaces, denoted by $\mathcal{M}_k(\Gamma_0(N),\chi)$ (resp. $\mathcal{S}_k(\Gamma_0(N),\chi)$). Moreover, $\mathcal{M}_k(\Gamma_0(N),\chi)$ is the direct sum of $\mathcal{S}_k(\Gamma_0(N),\chi)$ and Eisenstein series.

Every modular form $f(\tau)\in \mathcal{M}_k(\Gamma_0(N),\chi)$ admits a Fourier expansion at infinity of the form
$$f(\tau)=\sum_{n\ge 0}a(n)q^n$$
where $q:=e^{2\pi i\tau}$. Since spaces of modular forms are finite dimensional, given two modular forms with the same level $N$ and weight $k$, it is known that they are equal if their Fourier expansions agree for the first $k[SL_2(\mathbb{Z}):\Gamma_0(N)]/12$ terms (see \cite{Kob1993, Shi1971}).

For more details on modular forms, the reader may refer to \cite{DS2005}.

\section{Formulas of $A_{3,k}(n)$}\label{sec:3}

Taking $t=3$ in \eqref{eq:1.1}, we have
\begin{equation}\label{eq:3.1}
\Phi(q):=\sum_{n\ge 0}A_3(n)q^n=\frac{(q^3;q^3)_\infty^3}{(q;q)_\infty},
\end{equation}
which is essentially a quotient of Dedekind $\eta$-functions, where $\eta$ is defined by
$$\eta(\tau)=q^{1/24}\prod_{n\ge 1}(1-q^n)=q^{1/24}(q;q)_\infty$$
with $q:=e^{2\pi i\tau}$. In fact, we have
\begin{equation}\label{eq:3.2}
\Phi(q)=q^{-1/3}\frac{\eta^3(3\tau)}{\eta(\tau)}.
\end{equation}

Now we focus on $\Phi^k(q)$. To make our argument more complete, we include the proofs of \eqref{eq:wang1}, \eqref{eq:wang2}, and \eqref{eq:GO} as part of this section. We should mention at first that for odd $k$ our modular forms have nebentypus $\chi=(-3|n)$, whereas the nebentypus is $\chi=(9|n)$ when $k$ is even (here $(\cdot|n)$ denotes the Jacobi symbol).

\subsection{The case $\boldsymbol{k=1}$}

We consider the weight $1$ modular form $q\Phi(q^3)$ on $\Gamma_0(9)$ defined by
$$q\Phi(q^3)=\frac{\eta^3(9\tau)}{\eta(3\tau)}=\sum_{n\ge 0}A_{3,1}(n)q^{3n+1}.$$
The first few terms of its Fourier expansion are
$$q\Phi(q^3)=q+q^4+2q^7+2q^{13}+q^{16}+\cdots.$$
Let $\chi_3(n)=(n|3)$ denote the Legendre symbol, and write
$$\sigma_{0,\chi_3}(n)=\sum_{d\mid n}\chi_3(d).$$
It turns out that $q\Phi(q^3)$ is the Eisenstein series on $\Gamma_0(9)$ given by
$$q\Phi(q^3)=\sum_{n\ge 0}\sigma_{0,\chi_3}(3n+1)q^{3n+1}.$$
We therefore obtain
\begin{theorem}[Granville and Ono]
For any integer $n\ge 0$,
\begin{equation}\label{eq:3.1.main}
A_{3,1}(n)=\sigma_{0,\chi_3}(3n+1).
\end{equation}
\end{theorem}

\subsection{The case $\boldsymbol{k=2}$}

We consider the modular form $q^2 \Phi^2(q^3)\in \mathcal{M}_2(\Gamma_0(9))$.
$$q^2 \Phi^2(q^3)=\frac{\eta^6(9\tau)}{\eta^2(3\tau)}=\sum_{n\ge 0}A_{3,2}(n)q^{3n+2}.$$
The first few terms of its Fourier expansion are
$$q^2 \Phi^2(q^3)=q^2 + 2 q^5 + 5 q^8 + 4 q^{11} + 8 q^{14}+\cdots.$$
Let
$$\sigma_{1}(n)=\sum_{d\mid n}d.$$
We define the following weight $2$ Eisenstein series on $\Gamma_0(9)$
$$\sum_{n\ge 0}\sigma_{1}(3n+2)q^{3n+2}.$$
One readily verifies our generating function $q^2 \Phi^2(q^3)$ satisfies
$$q^2 \Phi^2(q^3)=\frac{1}{3}\sum_{n\ge 0}\sigma_{1}(3n+2)q^{3n+2}.$$
We therefore prove
\begin{theorem}[Wang]
For any integer $n\ge 0$,
\begin{equation}\label{eq:3.2.main}
A_{3,2}(n)=\frac{1}{3}\sigma_{1}(3n+2).
\end{equation}
\end{theorem}

\subsection{The case $\boldsymbol{k=3}$}

Here we consider the weight $3$ modular form $q\Phi^3(q)$ on $\Gamma_0(3)$:
$$q\Phi^3(q)=\frac{\eta^9(3\tau)}{\eta^3(\tau)}=\sum_{n\ge 0}A_{3,3}(n)q^{n+1}.$$
We give the first few terms of its Fourier expansion as follows
$$q\Phi^3(q)=q + 3 q^2 + 9 q^3 + 13 q^4 + 24 q^5 +\cdots.$$
Let
$$\sigma_{2,\chi_3}^{*}(n)=\sum_{d\mid n}\chi_3(n/d)d^2.$$
Now we consider the weight $3$ Eisenstein series on $\Gamma_0(3)$
$$\sum_{n\ge 1}\sigma_{2,\chi_3}^{*}(n)q^n.$$
By equating Fourier coefficients we find that
$$q\Phi^3(q)=\sum_{n\ge 1}\sigma_{2,\chi_3}^{*}(n)q^n.$$
It follows
\begin{theorem}[Wang]
For any integer $n\ge 0$,
\begin{equation}\label{eq:3.3.main}
A_{3,3}(n)=\sigma_{2,\chi_3}^{*}(n+1).
\end{equation}
\end{theorem}

\subsection{The case $\boldsymbol{k=4}$}

We consider the weight $4$ modular form $q^4\Phi^4(q^3)$ on $\Gamma_0(9)$ given by
$$q^4\Phi^4(q^3)=\frac{\eta^{12}(9\tau)}{\eta^4(3\tau)}=\sum_{n\ge 0}A_{3,4}(n)q^{3n+4}.$$
Here are the first few terms of the Fourier expansion of $q^4\Phi^4(q^3)$:
$$q^4\Phi^4(q^3)=q^4 + 4 q^7 + 14 q^{10} + 28 q^{13} + 57 q^{16}+\cdots.$$
Let
$$\sigma_3(n)=\sum_{d\mid n}d^3.$$
We consider the following weight $4$ Eisenstein series on $\Gamma_0(9)$:
$$E(\tau)=\sum_{n\ge 0}\sigma_3(3n+1)q^{3n+1}.$$
Note also that the space of cusp forms $\mathcal{S}_4(\Gamma_0(9))$ is $1$ dimensional and is spanned by the $\eta$-product
$$\eta^8(3\tau)=q - 8 q^4 + 20 q^7 - 70 q^{13} + 64 q^{16} +\cdots.$$
By equating Fourier coefficients we obtain the following identity:
$$q^4\Phi^4(q^3)=\frac{1}{81}(E(\tau)-\eta^8(3\tau)).$$
It implies
\begin{theorem}
Let $\eta^8(3\tau)=\sum_{n\ge 1}a(n)q^n$. For any integer $n\ge 0$,
\begin{equation}\label{eq:3.4.main}
A_{3,4}(n)=\frac{1}{81}(\sigma_3(3n+4)-a(3n+4)).
\end{equation}
\end{theorem}

\begin{remark}
It is known that $\eta^8(3\tau)$ is a cusp with complex multiplication. Since all forms with complex multiplication are lacunary, that is, the arithmetic density of their non-zero Fourier coefficients is $0$, we immediately see $A_{3,4}(n)=\sigma_3(3n+4)/81$ almost always.
\end{remark}

\subsection{The case $\boldsymbol{k=6}$}

Here we consider the modular form $q^2 \Phi^6(q)\in \mathcal{M}_6(\Gamma_0(3))$.
$$q^2 \Phi^6(q)=\frac{\eta^{18}(3\tau)}{\eta^6(\tau)}=\sum_{n\ge 0}A_{3,6}(n)q^{n+2}.$$
We give the first few terms of its Fourier expansion
$$q^2 \Phi^6(q)=q^2 + 6 q^3 + 27 q^4 + 80 q^5 + 207 q^6 +\cdots.$$
Given a prime $p$, let $\nu_p(n)$ denote the largest integer $e$ such that $p^e\mid n$. Write
$$\sigma_{5,3}^{\#}(n)=3^{5\nu_3(n)}\sum_{d\mid n\atop d\not\equiv 0\bmod{3}}d^5.$$
Now we define the following weight $6$ Eisenstein series on $\Gamma_0(3)$
$$E(\tau)=\sum_{n\ge 1}\sigma_{5,3}^{\#}(n)q^n.$$
It is easy to see $\mathcal{S}_6(\Gamma_0(3))$ is $1$ dimensional and is spanned by
$$\eta^6(\tau)\eta^6(3\tau)=q - 6 q^2 + 9 q^3 + 4 q^4 + 6 q^5+\cdots.$$
Our form $q^2 \Phi^6(q)$ satisfies
$$q^2 \Phi^6(q)=\frac{1}{39}(E(\tau)-\eta^6(\tau)\eta^6(3\tau)).$$
We thus conclude
\begin{theorem}
Let $\eta^6(\tau)\eta^6(3\tau)=\sum_{n\ge 1}a(n)q^n$. For any integer $n\ge 0$,
\begin{equation}\label{eq:3.6.main}
A_{3,6}(n)=\frac{1}{39}(\sigma_{5,3}^{\#}(n+2)-a(n+2)).
\end{equation}
\end{theorem}

\section{Formula of $A_{4,2}(n)$}

Taking $t=4$ in \eqref{eq:1.1}, we have
\begin{equation}\label{eq:4.1}
\Phi(q):=\sum_{n\ge 0}A_4(n)q^n=\frac{(q^4;q^4)_\infty^4}{(q;q)_\infty}.
\end{equation}
One readily sees it equals the following $\eta$-quotient
\begin{equation}\label{eq:4.2}
\Phi(q)=q^{-5/8}\frac{\eta^4(4\tau)}{\eta(\tau)}.
\end{equation}

Now we consider the modular form $q^5 \Phi^2(q^4)\in \mathcal{M}_3(\Gamma_0(16),(-4|n))$.
$$q^5 \Phi^2(q^4)=\frac{\eta^{8}(16\tau)}{\eta^2(4\tau)}=\sum_{n\ge 0}A_{4,2}(n)q^{4n+5}.$$
The first few terms of its Fourier expansion are
$$q^5 \Phi^2(q^4)=q^5 + 2 q^9 + 5 q^{13} + 10 q^{17} + 12 q^{21} +\cdots.$$
Let $\chi_{4,2}$ be the Dirichlet character mod $4$ given by
$$\chi_{4,2}(n)=\begin{cases}1 & n\equiv 1\bmod 4,\\-1 & n\equiv 3\bmod 4,\\0 & (n,4)> 1,\end{cases}$$
and write
$$\sigma_{2,\chi_{4,2}}(n)=\sum_{d\mid n}\chi_{4,2}(d)d^2.$$
We consider the following weight $3$ Eisenstein series on $\Gamma_0(16)$:
$$E(\tau)=\sum_{n\ge 0}\sigma_{2,\chi_{4,2}}(4n+1)q^{4n+1}.$$
It is also known that the space of cusp forms $\mathcal{S}_3(\Gamma_0(16),(-4|n))$ is $1$ dimensional and is spanned by
$$\eta^6(4\tau)=q - 6 q^5 + 9 q^9 + 10 q^{13} - 30 q^{17} +\cdots.$$
One readily verifies
$$q^5 \Phi^2(q^4)=\frac{1}{32}(E(\tau)-\eta^6(4\tau)).$$
We therefore prove
\begin{theorem}
Let $\eta^6(4\tau)=\sum_{n\ge 1}a(n)q^n$. For any integer $n\ge 0$,
\begin{equation}\label{eq:4.2.main}
A_{4,2}(n)=\frac{1}{32}(\sigma_{2,\chi_{4,2}}(4n+5)-a(4n+5)).
\end{equation}
\end{theorem}

\section{Formulas of $A_{5,k}(n)$}

Taking $t=5$ in \eqref{eq:1.1}, we have
\begin{equation}\label{eq:5.1}
\Phi(q):=\sum_{n\ge 0}A_5(n)q^n=\frac{(q^5;q^5)_\infty^5}{(q;q)_\infty}.
\end{equation}
It is easy to obtain the following identity:
\begin{equation}\label{eq:5.2}
\Phi(q)=q^{-1}\frac{\eta^5(5\tau)}{\eta(\tau)}.
\end{equation}

\subsection{The case $\boldsymbol{k=1}$}

Here we consider the weight $2$ modular form $q \Phi(q)$ on $\Gamma_0(5)$ with nebentypus $(5|n)$.
$$q \Phi(q)=\frac{\eta^{5}(5\tau)}{\eta(\tau)}=\sum_{n\ge 0}A_{5,1}(n)q^{n+1}.$$
The first few terms of its Fourier expansion are
$$q \Phi(q)=q + q^2 + 2 q^3 + 3 q^4 + 5 q^5 +\cdots.$$
Let $\chi_{5,3}$ be the Dirichlet character mod $5$ given by
$$\chi_{5,3}(n)=\begin{cases}1 & n\equiv 1,4\bmod 5,\\-1 & n\equiv 2,3\bmod 5,\\0 & (n,5)> 1.\end{cases}$$
Now write
$$\sigma_{1,\chi_{5,3}}^*(n)=\sum_{d\mid n}\chi_{5,3}(n/d)d.$$
Define the following weight $2$ Eisenstein series on $\Gamma_0(5)$
$$E(\tau)=\sum_{n\ge 1}\sigma_{1,\chi_{5,3}}^*(n)q^n.$$
It follows by equating Fourier coefficients that
$$q \Phi(q)=E(\tau).$$
We conclude
\begin{theorem}\label{th:5.1}
For any integer $n\ge 0$,
\begin{equation}\label{eq:5.1.main}
A_{5,1}(n)=\sigma_{1,\chi_{5,3}}^*(n+1).
\end{equation}
\end{theorem}

\begin{remark}
It is of interest to mention that we can prove Theorem \ref{th:5.1} through Bailey's $_6\psi_6$ formula. Recall
\begin{lemma}[Bailey's $_6\psi_6$ formula]
For $|qa^2/(bcde)| < 1$,
\begin{equation}\label{eq:Bailey}
\begin{split}
&_{6}\psi_{6}\left[\begin{matrix} q\sqrt{a},& -q \sqrt{a},& b,& c,& d, &e\\ \sqrt{a},& -\sqrt{a},& aq/b, &aq/c, &aq/d, &aq/e \end{matrix}; q, \frac{qa^2}{bcde}\right]\\
&\quad=\frac{(aq,aq/(bc),aq/(bd),aq/(be),aq/(cd),aq/(ce),aq/(de),q,q/a;q)_{\infty}}{(aq/b,aq/c,aq/d,aq/e,q/b,q/c,q/d,q/e,qa^2/(bcde);q)_{\infty}},
\end{split}
\end{equation}
where the $_s\psi_s$ function is defined as
$$_{s}\psi_{s}\left[\begin{matrix} a_1,\ldots,a_s\\ b_1,\ldots,b_s \end{matrix}; q, z\right]:=\sum_{n=-\infty}^{\infty}\frac{(a_1,\ldots,a_s;q)_n}{(b_1,\ldots,b_s;q)_n}z^n.$$
\end{lemma}
\noindent For its proof, the reader may refer to \cite[Sec. 5.3]{GR2004}. Now taking
$$(a,b,c,d,e,q)\to(q^4,q,q,q^3,q^3,q^5)$$
in \eqref{eq:Bailey}, we deduce that
\begin{align*}
&\sum_{n\ge 0}A_{5,1}(n)q^{n+1}=q\frac{(q^5;q^5)_\infty^5}{(q;q)_\infty}\\
&\quad=\sum_{n\ge 0}\left\{\frac{(1-q^{10n+4})(q^{5n+1}-q^{5n+3})}{(1-q^{5n+1})^2(1-q^{5n+3})^2}-\frac{(1-q^{10n+6})(q^{5n+2}-q^{5n+4})}{(1-q^{5n+2})^2(1-q^{5n+4})^2}\right\}\\
&\quad=\sum_{n\ge 0}\left\{\frac{q^{5n+1}}{(1-q^{5n+1})^2}-\frac{q^{5n+2}}{(1-q^{5n+2})^2}-\frac{q^{5n+3}}{(1-q^{5n+3})^2}+\frac{q^{5n+4}}{(1-q^{5n+4})^2}\right\}.
\end{align*}
Note that for $|q|<1$ we have
$$\frac{q}{(1-q)^2}=\sum_{d\ge 1}dq^d.$$
It therefore follows
$$\sum_{n\ge 0}A_{5,1}(n)q^{n+1}=\sum_{n\ge 0}\sum_{d\ge 1}\left\{dq^{(5n+1)d}-dq^{(5n+2)d}-dq^{(5n+3)d}+dq^{(5n+4)d}\right\}.$$
By equating coefficients we have
$$A_{5,1}(n)=\sum_{d\mid n+1}\chi_{5,3}(d)\frac{n+1}{d}=\sum_{d\mid n+1}\chi_{5,3}\left(\frac{n+1}{d}\right)d=\sigma_{1,\chi_{5,3}}^*(n+1).$$
\end{remark}

\subsection{The case $\boldsymbol{k=2}$}

We consider the modular form $q^2 \Phi^2(q)\in \mathcal{M}_4(\Gamma_0(5))$.
$$q^2 \Phi^2(q)=\frac{\eta^{10}(5\tau)}{\eta^2(\tau)}=\sum_{n\ge 0}A_{5,2}(n)q^{n+2}.$$
The first few terms of its Fourier expansion are
$$q^2 \Phi^2(q)=q^2 + 2 q^3 + 5 q^4 + 10 q^5 + 20 q^6 +\cdots.$$
Let
$$\sigma_{3,5}^{\#}(n)=5^{3\nu_5(n)}\sum_{d\mid n\atop d\not\equiv 0\bmod{5}}d^3,$$
where $\nu_5(n)$ denotes the largest integer $e$ such that $5^e\mid n$. Now we define the following weight $4$ Eisenstein series on $\Gamma_0(5)$
$$E(\tau)=\sum_{n\ge 1}\sigma_{3,5}^{\#}(n)q^n.$$
We also notice that the space of cusp forms $\mathcal{S}_4(\Gamma_0(5))$ is $1$ dimensional and is spanned by
$$\eta^4(\tau)\eta^4(5\tau)=q - 4 q^2 + 2 q^3 + 8 q^4 - 5 q^5+\cdots.$$
Our form $q^2 \Phi^2(q)$ satisfies
$$q^2 \Phi^2(q)=\frac{1}{13}(E(\tau)-\eta^4(\tau)\eta^4(5\tau)).$$
We thus conclude
\begin{theorem}
Let $\eta^4(\tau)\eta^4(5\tau)=\sum_{n\ge 1}a(n)q^n$. For any integer $n\ge 0$,
\begin{equation}\label{eq:5.2.main}
A_{5,2}(n)=\frac{1}{13}(\sigma_{3,5}^{\#}(n+2)-a(n+2)).
\end{equation}
\end{theorem}

\section{Formula of $A_{7,1}(n)$}

Taking $t=7$ in \eqref{eq:1.1}, we have
\begin{equation}\label{eq:7.1}
\Phi(q):=\sum_{n\ge 0}A_7(n)q^n=\frac{(q^7;q^7)_\infty^7}{(q;q)_\infty}.
\end{equation}
One easily finds
\begin{equation}\label{eq:7.2}
\Phi(q)=q^{-2}\frac{\eta^7(7\tau)}{\eta(\tau)}.
\end{equation}

Now we consider the modular form $q^2 \Phi(q)\in \mathcal{M}_3(\Gamma_0(7),(-7|n))$.
$$q^2 \Phi(q)=\frac{\eta^{7}(7\tau)}{\eta(\tau)}=\sum_{n\ge 0}A_{7,1}(n)q^{n+2}.$$
The first few terms of its Fourier expansion are
$$q^2 \Phi(q)=q^2 + q^3 + 2 q^4 + 3 q^5 + 5 q^6 +\cdots.$$
Let $\chi_{7,4}$ be the Dirichlet character mod $7$ given by
$$\chi_{7,4}(n)=\begin{cases}1 & n\equiv 1,2,4\bmod 7,\\-1 & n\equiv 3,5,6\bmod 7,\\0 & (n,7)> 1.\end{cases}$$
We write
$$\sigma_{2,\chi_{7,4}}^*(n)=\sum_{d\mid n}\chi_{7,4}(n/d)d^2.$$
Consider the following weight $3$ Eisenstein series on $\Gamma_0(7)$:
$$E(\tau)=\sum_{n\ge 1}\sigma_{2,\chi_{7,4}}^*(n)q^n.$$
Note also that $\mathcal{S}_3(\Gamma_0(7),(-7|n))$ is $1$ dimensional and is spanned by
$$\eta^3(\tau)\eta^3(7\tau)=q - 3 q^2 + 5 q^4 - 7 q^7 - 3 q^8 +\cdots.$$
It is easy to verify that
$$q^2 \Phi(q)=\frac{1}{8}(E(\tau)-\eta^3(\tau)\eta^3(7\tau)).$$
We therefore prove
\begin{theorem}
Let $\eta^3(\tau)\eta^3(7\tau)=\sum_{n\ge 1}a(n)q^n$. For any integer $n\ge 0$,
\begin{equation}\label{eq:7.1.main}
A_{7,1}(n)=\frac{1}{8}(\sigma_{2,\chi_{7,4}}^*(n+2)-a(n+2)).
\end{equation}
\end{theorem}

\subsection*{Acknowledgments}

The author thanks Wei Lin and Yucheng Liu for helpful discussions. The author also thanks the referee for valuable suggestions which have improved the readability of this paper.

\bibliographystyle{amsplain}

\begin{thebibliography}{99}

\bibitem{BN2014}
N. D. Baruah and K. Nath, Some results on $3$-cores, \textit{Proc. Amer. Math. Soc.} \textbf{142} (2014), no. 2, 441--448.

\bibitem{BN2015}
N. D. Baruah and K. Nath, Infinite families of arithmetic identities and congruences for bipartitions with $3$-cores, \textit{J. Number Theory} \textbf{149} (2015), 92--104. 

\bibitem{Boy2002}
M. Boylan, Congruences for $2^t$-core partition functions, \textit{J. Number Theory} \textbf{92} (2002), no. 1, 131--138.

\bibitem{Chen2013}
S. C. Chen, Congruences for $t$-core partition functions, \textit{J. Number Theory} \textbf{133} (2013), no. 12, 4036--4046.

\bibitem{Dai2015}
H. B. Dai, Arithmetic of $3^t$-core partition functions, \textit{Integers} \textbf{15} (2015), Paper No. A7, 5 pp.

\bibitem{DS2005}
F. Diamond and J. Shurman, \textit{A first course in modular forms}, Graduate Texts in Mathematics, \textbf{228}. Springer-Verlag, New York, 2005. xvi+436 pp.

\bibitem{GKS1990}
F. Garvan, D. Kim, and D. Stanton, Cranks and $t$-cores, \textit{Invent. Math.} \textbf{101} (1990), no. 1, 1--17.

\bibitem{GR2004}
G. Gasper and M. Rahman, \textit{Basic hypergeometric series. Second edition}, Encyclopedia of Mathematics and its Applications, \textbf{96}. Cambridge University Press, Cambridge, 2004. xxvi+428 pp. 

\bibitem{GO1996}
A. Granville and K. Ono, Defect zero $p$-blocks for finite simple groups, \textit{Trans. Amer. Math. Soc.} \textbf{348} (1996), no. 1, 331--347. 

\bibitem{HS1996}
M. D. Hirschhorn and J. A. Sellers, Some amazing facts about $4$-cores, \textit{J. Number Theory} \textbf{60} (1996), no. 1, 51--69. 

\bibitem{Kob1993}
N. Koblitz, \textit{Introduction to elliptic curves and modular forms. Second edition}, Graduate Texts in Mathematics, \textbf{97}. Springer-Verlag, New York, 1993. x+248 pp.

\bibitem{Lin2014}
B. L. S. Lin, Some results on bipartitions with $3$-core, \textit{J. Number Theory} \textbf{139} (2014), 44--52. 

\bibitem{ORW1995}
K. Ono, S. Robins, and P. T. Wahl, On the representation of integers as sums of triangular numbers, \textit{Aequationes Math.} \textbf{50} (1995), no. 1-2, 73--94. 

\bibitem{Shi1971}
G. Shimura, \textit{Introduction to the arithmetic theory of automorphic functions}, Reprint of the 1971 original. Publications of the Mathematical Society of Japan, \textbf{11}. Kan\^o Memorial Lectures, \textbf{1}. Princeton University Press, Princeton, NJ, 1994. xiv+271 pp.

\bibitem{Wang2016}
L. Wang, Explicit formulas for partition pairs and triples with $3$-cores, \textit{J. Math. Anal. Appl.} \textbf{434} (2016), no. 2, 1053--1064.

\bibitem{Yao2015}
O. X. M. Yao, Infinite families of congruences modulo $3$ and $9$ for bipartitions with $3$-cores, \textit{Bull. Aust. Math. Soc.} \textbf{91} (2015), no. 1, 47--52. 

\end{thebibliography}

\end{document}